\documentclass[a4ja,10pt,leqno]{article}
\usepackage{amsmath}
\usepackage{amsfonts}
\usepackage{latexsym}
\usepackage{amssymb}
\usepackage{amsthm}
\usepackage{amscd}
\newcommand{\sect}{\gamma}
\newcommand{\dbar}{d\!\!\lower-0.8ex\hbox{-}}

\newcommand{\integ}{{\scriptstyle\int}\kern-0.7em{\star}\,\,}

\newcommand{\symsum}{\hbox{$\sum
     $\kern-1em\lower-0.2ex\hbox{$\scriptscriptstyle{\bigcirc}$}}\,}
\newcommand{\ott}{\lower-0.4ex\hbox{${\scriptscriptstyle{\otimes}}$}}
\newcommand{\btt}{\lower-0.2ex\hbox{${\scriptscriptstyle{\bullet}}$}}
\newcommand{\dtt}{\lower-0.2ex\hbox{${\scriptscriptstyle{\odot}}$}}

\newtheorem{thm}{Theorem}[section]
\newtheorem{prop}[thm]{Proposition}

\newtheorem{lem}[thm]{Lemma}

\newtheorem{defn}[thm]{Definition}

\theoremstyle{remark}

\newcommand{\C}{\mathbb{C}}

\newcommand{\cstar}{{\mathbb C}_*[\zeta]}
\newcommand{\stk}{*_\kappa}

\newcommand{\p}{{\mathcal P}}

\title{Geometric objects in an approach to quantum geometry}

\author{Hideki Omori\\
Department of Mathematics, 
Tokyo University of Science, \\
Noda, Chiba, 278-8510, Japan, \\
\texttt{omori@ma.noda.tus.ac.jp}
\and
Yoshiaki Maeda
\footnote{Partially supported by Grant-in-Aid for 
Scientific Research (\#12440022.), Ministry of Education , Science and 
Culture, Japan.}\\
Department of Mathematics, 
Faculty of Science and Technology, \\
Keio University, Hiyoshi, Yokohama, 
223-8825, Japan, \\
\texttt{maeda@math.keio.ac.jp}; 
\and
Naoya Miyazaki
\footnote{
Partially supported by Grant-in-Aid for 
Scientific Research (\#15740045.), Ministry of Education , Science and 
Culture, Japan.}
\\
Department of Mathematics, Faculty of Economics, \\
Keio University, Hiyoshi, Yokohama, 
223-8521, Japan, \\
\texttt{miyazaki@math.hc.keio.ac.jp}; 
\and
Akira Yoshioka
\footnote{Partially supported by Grant-in-Aid for
  Scientific Research (\#11640095.), Ministry of Education , Science and 
  Culture, Japan.}
\\
Department of Mathematics, 
Tokyo University of Science,\\
Kagurazaka, Tokyo, 102-8601, Japan,\\
\texttt{yoshioka@rs.kagu.tus.ac.jp}; 
}

\begin{document}
\maketitle

\begin{center}
\mbox{\parbox{.8\linewidth}{
{\bf Abstract} 
Ideas from deformation quantization applied to algebras with one
generator lead to methods to treat a nonlinear flat connection. 
It provides us elements of algebras to be parallel sections. The
moduli space of the parallel sections is studied as an example of
bundle-like objects with discordant (sogo) transition functions,    
which suggests us to treat movable branching singularities.   
 }}  
\end{center}
\par\noindent
{\small\bf A.M.S Classification (2000)}: 
{Primary 53D55, 53D10; Secondary 46L65}

\section{Introduction}
The aim of this paper is to show that deformation quantization 
provides us a new geometric idea going beyond classical geometry.  
In fact, there has been several 
attempts to describe ``quantum objects'' in a geometric way
(cf. \cite{Co}, \cite{landi}, \cite{manin}), although no treatment has
been accepted as definitions.  
Motivated by these attempts, we produce a proposal to describe
objects which arise from the deformation of algebras, 
as one of approaches to describe quantum mechanics mathematically is  
the notion of deformation quantization, which was a deformation of 
Poisson algebras. 
Through the construction of the star 
exponential functions of the quadratic forms in the complex Weyl 
algebra, we found several strange phenomena, which cannot be treated  
as classical geometric objects 
(cf.\,\cite{om1}, \cite{ommy}, \cite{OMMY5}, \cite{OMMY6}).  
Our main concern is to answer how to handle these objects  
geometrically, and we consider our results are a step toward quantum geometry. 
However, similar questions arise even for deformations of commutative
algebras, as for the case deformation quantizations.
For this reason in this paper, we deal with the simplest case of the 
deformation of the associative commutative algebra of polynomials of one variable.  

\smallskip
Let $\C_{*}[\zeta ]$ be the abstract algebra over ${\C}$ 
generated by $\zeta$, which is trivially a commutative associative algebra.  
$\C_*[\zeta]$ is canonically isomorphic to the ordinary commutative 
algebra ${\C}[\zeta ]$ of polynomials of $\zeta$. 

Motivated by deformation quantization, 
we introduce associative commutative products on ${\C}[\zeta ]$ 
parametrized by a complex number $\kappa$ (cf.\eqref{star-product}), 
which gives both a deformation of the canonical product and a 
representation with the parameter $\kappa$ of $\C$. 
 
Our standpoint formulated in \S\,\ref{1.1} is to view elements in the abstract algebra 
${\mathbb C}_* [\zeta]$ as a family of elements. The deformation parameter $\kappa$ is viewed 
as an indeterminate. 

One method of treating this family of elements as geometric objects is to
introduce the notion of infinitesimal intertwiners, which play the  
role of a connection. In fact, elements of ${\mathbb C}_* [\zeta]$  
can be viewed as parallel sections with respect to this connection. 
These elements elements are called $q$-number polynomials. 

In \S\,\ref{1.2} and \S\,\ref{1.3}, we extend this settings to a class 
of transcendental elements such as exponential functions. 
In this setting, the notion of densely defined 
multi-valued parallel sections appears crucially. 
We also call these $q$-number functions in analogy with  \cite{AW}.  
However, the only geometrical setting is to extend the infinitesimal 
intertwiners to a linear connection on a trivial 
bundle over $\C$ with a certain Fr{\'e}chet space of entire functions.   

\smallskip
In \S\,\ref{3.0} we investigate the moduli space of densely 
defined parallel sections consisting of exponential functions of
quadratic forms. We show that the moduli space is not an ordinary bundle, 
as it contains fuzzy transition functions. This has similarities to
the theory of gerbes (cf. \cite{brylinski}, \cite{murr}). 

However, our construction has a different flavor from the differential 
geometric point of view, since gerbes are classified by 
the Dixmir-Douady class in the third cohomology, while our example is
constructed on the 2-sphere or the complex plane. We would like to 
call this fuzzy object a {\it pile}, which is roughly a collection of
a local trivializations without the consistency property for the transition 
functions.  

\smallskip

We come up a similar situation in quantizing non-integrable 
closed 2-forms on manifolds. 
As for integrable symplectic forms on symplectic manifolds, 
we can construct a prequantum bundle, which is a line
bundle with the connection whose curvature is given by
the symplectic form. 

We attempt the prequantization of a non-integrable closed 2-form 
by stimulating our examples describing the moduli space of densely 
defined multi-valued parallel sections.
We note that Melrose \cite{melrose} proposed a method handling a type of 
prequantization of non-integrable closed 2-forms, 
which seems closely related with our approach. 

We note that the deformation parameter $\kappa$ can be 
viewed as an indeterminate moving in $\C$. Parallel sections in our 
example have branching singular points with respect to $\kappa$,  
depending on the parallel section. Such movable branching singularities 
have never treated systematically in the theory of differential equations. 
However, when we treat $\kappa$ as an indeterminate, 
this movable singularity can be managed by algebraic calculations.     
The aim of the last section is to give simple examples.

We introduce an associative product on the space of parallel sections
of exponential functions of quadratic forms, but this product is 
``broken'' in the sense that for every $\kappa$, there is a singular set on 
which the product diverges. Thanks to the movable singularities, this 
broken product defines an associative product by treating $\kappa$ as 
an indeterminate. This computation provides a novel aspect of the 
noncommutative calculus.  

In the end, our work seems to extend the notion of {\it points} as 
established elements of a fixed set to a more flexible notion of elememts.

\section{Deformation of a commutative product}

Let $\cstar$ be the algebra over $\C$ generated by $\zeta$. 
For convenience, we denote by $*$ the product on the algebra $\cstar$. 
The algebra $\cstar$ is isomorphic to the algebra $\C[\zeta]$ of polynomials in 
$\zeta$ over $\C$, but we will view $\cstar$ as a family of several 
algebras which are mutually {\it isomorphic}.  

\subsection{A deformation of commutative product on ${\mathbb C}[\zeta]$}
\label{1.1}

We denote the set of polynomials of $\zeta$ by $\p({\mathbb C})$ when  
it is viewed as a linear space. We introduce a family of product $*_\kappa$ on 
$\p({\mathbb C})$ parametrized by $\kappa\in \C$ as follows. 

\begin{defn}
For every $f, g\in \p({\mathbb C})$, we set 
\begin{eqnarray}
\nonumber
f\stk g
=\sum_{\ell=0}^\infty\frac{1}{\ell!}\bigr(\frac{\kappa}{2}\bigr)^\ell 
\partial_\zeta^\ell f(\zeta) \cdot \partial_\zeta^\ell g (\zeta) . 
\label{star-product}
\end{eqnarray}
\end{defn}
Then, $(\p({\mathbb C}),*_\kappa)$ is an associative commutative 
algebra for every $\kappa\in {\mathbb C}$. 
Since putting $\kappa=0$ gives the algebra ${\mathbb C}[\zeta]$,  
the family of algebras 
$\{\p({\mathbb C}),*_\kappa)\}_{\kappa\in {\mathbb C}}$ 
gives a deformation of ${\mathbb C}[\zeta]$ within 
an associative commutative algebra. 
We note the following. 
\begin{lem}
For every $\kappa,\kappa'\in \C$, the algebras $(\p({\mathbb C}),*_\kappa)$ 
and $(\p({\mathbb C}),*_{\kappa'})$ 
are mutually isomorphic. 
Namely, the mapping 
$T_\kappa^{\kappa^\prime}:\p({\mathbb C})\rightarrow
\p({\mathbb C})$ given by  
\begin{equation}
\label{intertwiner}
 T_\kappa^{\kappa^\prime} (f)
 =\Big(\exp\frac{1}{4}(\kappa'{-}\kappa)\partial_\zeta^2\Big)f(\zeta)
=\sum_{\ell=0}^\infty\frac{1}{\ell!}\bigr(
\frac{1}{4}(\kappa'{-}\kappa)\bigr)^\ell(\partial_\zeta^{2\ell})f(\zeta)
\end{equation} 
satisfies 
$T_\kappa^{\kappa^\prime} (f*_\kappa g) = 
   T_\kappa^{\kappa^\prime} (f)*_{\kappa^\prime}
   T_\kappa^{\kappa^\prime} (g).$
 \end{lem}

\begin{defn}
The isomorphism $T_\kappa^{\kappa'}$ given by \eqref{intertwiner} 
is called the {\it intertwiner} 
between the algebras 
$(\p({\mathbb C}),*_\kappa)$ 
and $(\p({\mathbb C}),*_{\kappa'})$.
\end{defn}
Taking the derivative in $\kappa'$ for $T_\kappa^{\kappa'}$ defines an 
{\it infinitesimal intertwiner}.  
Namely we set for $\kappa \in {\mathbb C}$, 
\begin{eqnarray}
\label{infinitesimal}
t_\kappa(u)(f)=\frac{d}{ds}\bigr|_{s=0}T_{\kappa}^{\kappa+su}(f)=
\frac{1}{4}u\partial_\zeta^2 f.
\end{eqnarray}
The infinitesimal intertwiner gives a realization of $\cstar $ as follows. 
Let $\pi:{\mathbb C}{\times}\p({\mathbb C})\rightarrow {\mathbb C}$ be 
the trivial bundle over ${\mathbb C}$,  
and $\Gamma({\mathbb C}{\times}P({\mathbb C}))$   
the set of sections of this bundle.  
Using the infinitesimal intertwiner defined by 
\eqref{infinitesimal}, we introduce a connection $\nabla$ on 
$\Gamma({\mathbb C}{\times}\p({\mathbb C}))$: 
For a smooth curve $c(s)$ in $\C$ and $\sect \in 
\Gamma({\mathbb C}{\times}\p({\mathbb C}))$, 
we set 
\begin{equation}\label{def-connection}
\nabla_{\dot{c}}\sect(s){=}
\frac{d}{ds}\sect(c(s)){-}t_{c(s)}(\dot{c}(s))(\sect),
\quad {\text{where}}\,\,\,
\dot{c}(s){=}\frac{d}{ds}c(s).
\end{equation}

kkk
\begin{defn}\label{defn}
A section $\sect{\in}\Gamma({\mathbb C}{\times}\p({\mathbb C}))$ is 
{\it parallel} if $\nabla\sect{=}0$. 
We denote by ${\mathcal S}({\mathbb C}{\times}\p({\mathbb C}))$ the
set of all parallel sections 
$\sect{\in}\Gamma({\mathbb C}{\times}\p({\mathbb C}))$. 
\end{defn}
Let us consider an element $f_*{\in} {\mathbb C}_*[\zeta]$.
According to the unique expression of an element 
$f_*\in {\mathbb C}_*[\zeta]$ as 
\begin{equation*}
f_*=\sum a_j\underbrace{\zeta{*}\cdots{*}\zeta}_{j\mbox{-times}}
\quad \mbox{(finite sum)}, a_j{\in}{\mathbb C}, 
\end{equation*}
we set the element $f_{\kappa}\in \p({\mathbb C})$ for every 
$\kappa\in {\mathbb C}$ by 
\begin{equation*}
f_{\kappa}{=}
\sum a_j\underbrace{\zeta *_\kappa\cdots *_\kappa\zeta}_{j\mbox{-times}} 
\quad \mbox{(finite sum)}, a_j{\in}{\mathbb C}.
\end{equation*}
The section $\sect_{f_*}(\kappa){=}f_{\kappa}$ gives a parallel 
section of the bundle 
$\pi: {\mathbb C}{\times}\p({\mathbb C}){\rightarrow}{\mathbb C}$. 

Using the product formula $*_\kappa$, we define a product $*$ on 
${\mathcal S}({\mathbb C}{\times}\p({\mathbb C}))$ as 
\begin{equation}
\label{product-star}
(\sect_1{*}\sect_2)(\kappa){=}\sect_1(\kappa)*_\kappa\sect_2(\kappa),
\quad \sect_1,\,\,\sect_2\in 
{\mathcal S}
({\mathbb C}{\times}\p({\mathbb C})). 
\end{equation}
\begin{lem}
\label{}
$({\mathcal S}
({\mathbb C}{\times}\p({\mathbb C})),*)$ is an associative commutative algebra. 
\end{lem}
This procedure gives an identification of the algebra   
$({\mathcal S}({\mathbb C}{\times}\p({\mathbb C})),*)$ with $\cstar $. 
Elements of $({\mathcal S}({\mathbb C}{\times}\p({\mathbb C})),*)$ will
be called $q$-number polynomials. 

\subsection{Strange exponential functions}\label{1.2}
We now extend this procedure to exponential functions. 

For $f_* \in \C_*[\zeta]$, we want to describe the star exponential 
functions $\exp_*f_*$, which may be a very transcendental element. 

\smallskip
Let ${\mathcal E}({\C})$ be the set of all entire functions on $\C$. 
For $p>0$, we set 
\begin{eqnarray}
{\mathcal E}_{p}({\C})&=&\{f\in{\mathcal E}(\C)~|~||f||_{p,\delta}=
\sup_{\zeta\in \C}e^{-\delta|\zeta|^{p}}|f(\zeta)|<\infty,
\forall \delta>0\},
\label{entire-p}
\end{eqnarray}
and also set 
${\mathcal E}_{p+}(\C){=}\cap_{q>p}{\mathcal E}_q(\C).$ 
It is seen that $({\mathcal E}_{p}({\C}), *_\kappa)$  
is a Fr\'{e}chet commutative associative algebra for $p\le 2$ (cf. \cite{ommy}).   
Recalling the intertwiner $T_\kappa^{\kappa'}$ given by 
\eqref{intertwiner}, we have the following \cite{OMMY5}:  
\begin{lem}
Let $p\le 2$. 
The intertwiner $T_\kappa^{\kappa'}$ 
in \eqref{intertwiner} is canonically extended to the map 
$T_\kappa^{\kappa^\prime} 
: {\mathcal E}_p({\mathbb C}) \rightarrow {\mathcal E}_p({\mathbb C})$ 
satisfying 
\begin{equation}
\label{intertwiner-condition-E}
T_\kappa^{\kappa^\prime} (f*_\kappa g) {=} 
   T_\kappa^{\kappa^\prime}(f) *_{\kappa^\prime}
   T_\kappa^{\kappa^\prime}(g) \quad\text{for every} 
   \,\,f,\,\,g\in{\mathcal E}_p({\mathbb C}). 
\end{equation} 
\end{lem}
We note that while the product $*_\kappa$ 
and the intertwiner $T_\kappa^{\kappa'}$ 
does not give an associative commutative product 
and does not extend to ${\mathcal E}_p({\mathbb C})$ for $p\ge 2$, the notion 
of connection $\nabla$ is still defined. 

Namely, we consider the trivial bundle  
$\pi:{\mathbb C}{\times}{\mathcal E}({\mathbb C})\rightarrow{\mathbb C}$  
over ${\mathbb C}$ with the fiber ${\mathcal E}({\mathbb C})$,
and the set of sections $\Gamma({\mathbb C}{\times}{\mathcal E}({\mathbb C}))$. 
We define for 
$\sect{\in}\Gamma({\mathbb C}{\times}{\mathcal E}({\mathbb C}))$ 
a covariant derivative $\nabla_{\dot{c}}\sect$ 
as the natural extension of \eqref{def-connection}. 
It is easily seen that $\nabla$ is well defined for 
$\Gamma({\mathbb C}{\times}{\mathcal E}_p({\mathbb C}))$ and 
$\Gamma({\mathbb C}{\times}{\mathcal E}_{p+}({\mathbb C}))$ for 
every $p{\ge}0$.  
Similarly as before we denote 
by ${\mathcal S}({\mathbb C}{\times}{\mathcal E}_p({\mathbb C}))$, 
${\mathcal S}({\mathbb C}{\times}{\mathcal E}_{p+}({\mathbb C}))$ 
the set of parallel 
sections. 

We wish to treat the star exponential function $\exp_*f_*$ 
for $f_*\in {\mathbb C}_*[\zeta]$. 
As in the argument in \S2.1, we have the realization 
$\{f_\kappa\}_{\kappa\in {\mathbb C}}$ 
of $f_*{\in}\C_*[\zeta]$, where $f_\kappa{\in}\p({\mathbb C})$. 
Fixing the $*_\kappa$ product gives the star exponential functions 
of $f_\kappa\in \p({\mathbb C})$ with 
respect to $*_\kappa$. We consider the following evolution equation 
\begin{equation}
\label{evolution-equation}
\left\{\begin{array}{l}
\partial_tF_{\kappa}(t)=f_\kappa(\zeta)*_\kappa F_\kappa(t),\\
F_\kappa(0)=g_\kappa . 
\end{array}\right. 
\end{equation}
If the equation (\ref{evolution-equation}) has  
a real analytic solution in $t$, then this solution is unique. 
Thus, we may set $\exp_{*_\kappa}f_\kappa=F_\kappa (1)$ when
(\ref{evolution-equation}) has an analytic solution with $F_{\kappa}(0)=1$. 

By moving $\kappa{\in}\C$, the totality of the star exponential 
functions $\{\exp_{*_\kappa}f_\kappa\}_{\kappa\in{\mathbb C}}$ may be
viewed as a natural representation of the star exponential function $\exp_*f_*$.  

\medskip
As an example, we consider the linear function $f(\zeta)=a\zeta$, where 
$a\in \C$. Then, the evolution equation $(\ref{evolution-equation})$ is 
expressed as 
\begin{equation}
\label{linear-evolution-equation}
\left\{\begin{array}{l}
\partial_tF_{\kappa}(t)=a\zeta F_\kappa+\frac{\kappa}{2}a\partial_\zeta F_\kappa , \\
F_\kappa(0)=1  . 
\end{array}\right. 
\end{equation}
By a direct computation, we have 
\begin{lem}\label{lem-linear-evolution-solution}
The equation \eqref{linear-evolution-equation} has the solution 
$F_\kappa(t){=}\exp (at\zeta{+}\frac{\kappa}{4}a^2t^2)$.
Thus, we may set 
\begin{equation}
\label{def-exp-linear-evolution-solution}
\exp_{*_\kappa}t\zeta=\exp (t\zeta+\frac{\kappa}{4}t^2)
\end{equation}
which is contained in ${\mathcal E}_{1+}(\C)$ for every $\kappa\in\C$. 
\end{lem}
Since the intertwiner $T_\kappa^{\kappa'}$ is defined 
on ${\mathcal E}_{1+}({\mathbb C})$, and 
$T_\kappa^{\kappa'}(\exp_{*_\kappa}a\zeta)=\exp_{*_{\kappa'}}a\zeta$, 
we see that $\{\exp_{*_\kappa}a\zeta\}_{\kappa\in{\mathbb C}}$ 
is an element of 
${\mathcal S}({\mathbb C}{\times}{\mathcal E}_{1+}({\mathbb C}))$.  
As in the argument as in \S2.1, it is natural to regard 
$\{\exp_{*_\kappa} a \zeta\}_{\kappa\in{\mathbb C}}$ as 
the star exponential function $\exp_{*}a\zeta$, which may be called a
$q$-number exponential function. 

\medskip
From the star exponential functions $\exp_{*_\kappa}a\zeta$, 
we construct a type of {\it delta function} via the star Fourier transform: 
Namely, we call  
\begin{equation}\label{def-delta}
\delta_{*_\kappa}(\zeta)=\int_{-\infty}^\infty\exp_{*_\kappa}it\zeta\,dt
\end{equation}
the $*_\kappa$-delta function. Using 
(\ref{def-exp-linear-evolution-solution}), we have 
\begin{lem}
The $*_\kappa$-delta function $\delta_{*_\kappa}(\zeta)$ is well defined as an 
element of ${\mathcal E}_{2+}(\C)$ for every $\kappa{\in}\C$ such that  
${\rm Re}(\kappa)>0$.  
\end{lem}
Using integration by parts, we easily see that   
$$
e^{i\theta}\int_{-\infty}^{\infty}exp_{*_\kappa}{e^{i\theta}it\zeta}\,dt, \quad
\,\,{\rm{Re}}\,e^{2i\theta}\kappa{>}0
$$ 
does not depend on $\theta$ whenever ${\rm Re}(e^{2i\theta}\kappa)>0$.
This allows us to define $\delta_{*_\kappa}(\zeta){\in}{\mathcal E}_{2+}(\C)$
for  $\kappa{\in}\C{-}\{0\}$.

\begin{lem}\label{multi-valued}
The mapping $\delta_*:\C{-}\{0\}\rightarrow{\mathcal E}_{2+}(\C)$
defined by $\kappa{\to}\delta_{*_\kappa}$ is two-valued. 
\end{lem}
\begin{proof}
We set
\begin{equation}\label{kappa-delta}
\delta(\zeta; e^{i\theta}, \kappa)=
e^{i\theta}\int_{-\infty}^\infty
\exp(ie^{i\theta}t\zeta-\frac{\kappa}{4}e^{2i\theta}t^2)dt. 
\end{equation}

\noindent
\mbox{\parbox{.4\linewidth}{
\setlength{\unitlength}{.25mm}
\begin{picture}(200,150)
\thinlines
\put(0,60){\line(1,0){200}}
\put(100,0){\line(0,1){130}}
\put(100,50){\line(2,-1){100}}
\put(100,50){\line(-2,1){100}}
\put(100,70){\line(2,-1){100}}
\put(100,70){\line(-2,1){100}}
\put(100,130){$\theta$}
\put(190,60){$\tau$}
\put(105,58){$\bullet$}
\put(105,48){$\bullet$}
\thicklines
\put(100,60){\line(1,0){11}}
\put(110,60){\line(0,-1){10}}
\put(110,50){\line(1,0){18}}
\put(128,50){\line(0,-1){11}}
\put(128,40){\line(1,0){18}}
\put(146,40){\line(0,-1){11}}
\put(146,30){\line(1,0){18}}
\put(0,5){\footnotesize{Fig.1}}
\end{picture}
}}
\hfill{\parbox{.52\linewidth}{
Setting $\kappa{=}e^{i\tau}$ gives  
that $\delta(\zeta; e^{i\theta}, e^{i\tau})$ is well defined on
the strip $\bullet$ given in Fig.1. 
Note that $\delta(\zeta; e^{i\theta},\kappa)$ depends only on 
$\tau$ in the strip 
${-}{\pi}{2}{<}\tau{+}2\theta{<}{\pi}{2}$ and 
$\delta(\zeta; e^{i\theta},\kappa)$ is a parallel section with
respect to $\kappa$. 

By moving $\tau$ from $0$ to $2\pi$ and also $\theta$ 
such that $(\tau, \theta)$ is contained in this strip.
The continuous chasing from $\tau{=}0$ to $\tau{=}2\pi$ in the 
defining strip gives}}
$$
\delta(\zeta; 1,c){=}
\int_{-\infty}^{\infty}\exp(it\zeta{-}\frac{1}{4}ct^2)dt 
{=}
\int_{-\infty}^{\infty}\exp(-it\zeta{-}\frac{1}{4}ct^2)dt 
{=}{-}\delta(\zeta; 1,c).
$$
\end{proof}
Let us consider the trivial vector bundle over $\C{-}\{0\}$. 
Lemma~\ref{multi-valued} tells us that $\delta_*(\zeta)$ can be viewed 
as a double-valued holomorphic parallel section over $\C{-}\{0\}$. Note
that $\delta(\zeta; 1,c){=}\frac{1}{\sqrt{2\pi c}}e^{-\frac{1}{c}\zeta^2}$, 
and $\lim_{c\to 0}\delta(\zeta; 1,c)$ gives us the ordinary delta function.

As seen in the star delta function as above, we need 
the notion of {\it densely defined multi-valued parallel 
sections} when we wish to extend the notion of parallel section 
beyond the algebra 
$({\mathcal E}_2({\C}),*_\kappa)$.  

\subsection{Star exponential functions of the quadratic functions}\label{1.3}

We set  
$$
P^{(2)}({\mathbb C}) 
{=}\{f(\zeta){=}a\zeta^2{+}b\,|\,a, b {\in}{\mathbb C}\},\quad 
{\mathbb C}_*^{(2)}[{\mathbb C}]{=}\{f_*(\zeta){=}a\zeta{*}\zeta
\in\C_*[\zeta]  \,| \, a{\in}{\mathbb C}\}
$$
Thus, we view $a\zeta{*}\zeta$ as the section 
$\sect(\kappa)=a\zeta^2{+}\frac{a}{2}
\kappa{\in}\Gamma({\mathbb C}{\times}P^{(2)}({\mathbb C}))$, 
where 
$\pi:{\mathbb C}{\times}P^{(2)}({\mathbb C}){\rightarrow}{\mathbb C}$ 
is the trivial bundle over 
${\mathbb C}$ with the fiber $P^{(2)}({\mathbb C})$. 
We now attempt to give a a meaning of the star exponential function 
$\exp_*a\zeta{*}\zeta$, $a{\in}{\mathbb C}$ along the argument in 
\S~\ref{1.1}. 

We consider a quadratic element $f_{*}\in \C_*^{(2)}[\zeta]$.  
Then the corresponding polynomial $f_{\kappa}$ is given by  
\begin{equation}\label{quadratic-element}
f_\kappa=\zeta{*_\kappa}\zeta=\zeta^2{+}\frac{\kappa}{2}. 
\end{equation} 
As in \S~\ref{1.1}, we view $\{f_\kappa\}_{\kappa\in\C}$ as 
a parallel section of $\C{\times}P^{(2)} (\C)$.  
We consider the following evolution equation. 
\begin{equation}
\label{quad-evolution-equation}
\partial_tF_{\kappa}(t)=f_\kappa(\zeta)*_\kappa F_\kappa(t), \quad F_\kappa(0)=g_{\kappa}
\end{equation}
where $f_\kappa$ is given by (\ref{quadratic-element}). 
\eqref{quad-evolution-equation} is rewritten as  
\begin{equation}
\label{quadratic-evolution-equation}
\partial_tF_{\kappa}(t)
=(\zeta^2{+}\frac{\kappa}{2})F_\kappa{+}\kappa\zeta\partial_\zeta F_\kappa{+}
\frac{\kappa^2}{4}\partial_\zeta^2 F_\kappa,\quad F_\kappa(0)=g_{\kappa}.
\end{equation}
We assume that the initial condition $g_\kappa$ is given by the form 
$g_\kappa=\rho_{\kappa,0}\exp a_{\kappa,0}\zeta^2,$ 
where 
$\rho_{\kappa,0}\in {\mathbb C}_\times{=}\C{-}\{0\}$ and 
$a_{\kappa,0}\in {\mathbb C}$. 
Putting $g_\kappa{=}1$ gives the star exponential function 
$\exp_{*_\kappa}f_\kappa(\zeta)$. 
To solve (\ref{quadratic-evolution-equation}) explicitly, we assume that 
$F_\kappa$ is the following form: 
\begin{equation}
\label{assumption}
F_\kappa(t)=\rho_\kappa(t)\exp a_\kappa(t)\zeta^2.
\end{equation}
Plugging (\ref{assumption}) to (\ref{quadratic-evolution-equation}), 
we have 
\begin{equation}\label{quadratic-evolution-ordinary-equation}
\left\{\begin{array}{l}
\partial_ta_\kappa{=}1{+}2a_\kappa \kappa{+}a_\kappa^2\kappa^2 , \\
\partial_t\rho_\kappa{=}    
\frac{\kappa }{2}\bigr(1{+}\kappa a_\kappa\bigr)\rho, \\
a_\kappa(0){=}a_{\kappa,0},\quad\rho_\kappa(0){=}\rho_{\kappa,0}. 
\end{array}\right.
\end{equation}
\begin{prop}\label{quadratic-evolution-solution}
The solution of \eqref{quadratic-evolution-ordinary-equation} is given by 
\begin{equation}
\label{solution-a-kappa22}
a_\kappa(t)=
\frac{a_{\kappa ,0}{+}t(1{+} \kappa a_{\kappa,0})}
{1{-}\kappa t(1{+}\kappa a_{\kappa,0})}, 
\quad
\rho_\kappa(t)=
\frac{\rho_{\kappa,0}}
{\sqrt{1{-}\kappa t(1{+}\kappa a_{\kappa,0})}}, 
\end{equation}
\end{prop}

\noindent
where we mind an ambiguity for choosing 
the sign of the square root in \eqref{solution-a-kappa22}.  
We set the subset ${\mathcal E}^{(2)}({\mathbb C})$ of 
${\mathcal E}({\mathbb C})$ as 
\[{\mathcal E}^{(2)}({\mathbb C})=
\{f=\rho\exp a \zeta^2 \,|\, \rho\in {\mathbb C}_\times,\,a\in{\mathbb C}\}.
\]
Assigning $f{=}\rho\exp a\zeta^2\in{\mathcal E}^{(2)}({\mathbb C})$ with 
$(\rho,a)$ gives an identification ${\mathcal E}^{(2)}({\mathbb C}) 
\cong{\C}_{\times}{\times}{\C}$. 
Note that 
${\mathcal E}^{(2)}({\mathbb C})$ is not contained in 
${\mathcal E}_2({\mathbb C})$ but in ${\mathcal E}_{2+}({\mathbb C})$, 
on which the product $*_\kappa$ may give strange phenomena (cf.\cite{OMMY5}). 

Consider the trivial bundle 
$\pi:{\C}{\times}{\mathcal E}^{(2)}({\C})$ 
over ${\C}$ with the fiber ${\mathcal E}^{(2)}({\mathbb C})$. 
In particular, putting $a_{\kappa,0}{=}0$, $\rho_{\kappa,0}{=}1$ 
and $t{=}a$ in Proposition~\ref{quadratic-evolution-solution}, we see 
\begin{equation}
\label{might-set}
\exp_{*_\kappa}a\zeta{*_\kappa}\zeta =\frac{1}{\sqrt{1{-}a\kappa}}\exp 
\frac{a}{1{-}a\kappa}\zeta^2
\end{equation}
where the right hand side of 
(\ref{might-set}) 
still contains an ambiguity for choosing the sign of the square root. 

Keeping this ambiguity, we have a kind of fuzzy one 
parameter group property for the exponential function of 
(\ref{might-set}). 
Namely, we set $g_\kappa=\exp_{*_\kappa}b\zeta{*_\kappa}\zeta$, 
where $b\in{\mathbb C}$, the solutions 
of (\ref{quadratic-evolution-equation}) yields 
the exponential law: 
\begin{equation}\label{exponentiallaw}
        \exp_{*_\kappa} a\zeta{*_{\kappa}}\zeta
                *_\kappa 
                        \exp_{*_\kappa} b\zeta{*_{\kappa}}\zeta
        = \frac{1}{\sqrt{1{-}(a{+}b)\kappa}}
                e^{\frac{a{+}b}{1{-}(a{+}b)\kappa}\zeta^2}=
                   \exp_{*_\kappa}(a{+}b)\zeta{*_{\kappa}}\zeta,
\end{equation}
where (\ref{exponentiallaw}) still contains an ambiguity for choosing 
the sign of the square root. 

\medskip
Recall the connection $\nabla$ on the trivial bundle 
$\pi:{\mathbb C}{\times}{\mathcal E}({\mathbb C})
{\rightarrow}{\mathbb C}$. 
It is easily seen that the connection $\nabla$ parallelizes the bundle 
$\pi:{\mathbb C}{\times}{\mathcal E}^{(2)}({\mathbb C}){\rightarrow} 
{\mathbb C}$. 
According to the identification 
${\mathcal E}^{(2)}({\mathbb C})\cong{\C}_{\times}{\times}{\C}$, 
we write $\sect(\kappa)=\rho(\kappa)\exp a(\kappa)\zeta^2$ as 
$(\rho(\kappa), a(\kappa))$. Then the equation $\nabla_{\partial_t}\sect=0$ gives 
\begin{equation}\label{parallel-equation}
\left\{\begin{array}{ll}
\partial_ta(t)=a(t)^2, \\
\partial_t\rho(t)=-\frac{1}{2}\rho(t)a(t). 
\end{array}\right. \end{equation}
We easily see that (\ref{might-set}) gives a densely defined 
parallel section. As seen in \cite{OMMY5}, it should also be considered 
as a densely defined multi-valued section of this bundle.
Thus, we may view the star exponential function 
$\exp_*a\zeta{*}\zeta$ as a family 
$$
\big\{F_\kappa(\zeta)=\frac{1}{\sqrt{1{-}a\kappa}}
\exp\frac{a}{1{-}a\kappa}\zeta^2\big\}_{\kappa\in{\mathbb C}}.
$$

The realization of $\exp_*a\zeta{*}\zeta$ is a parallel section 
$\sect(\kappa)=\rho_\kappa\exp a\kappa\zeta^2$ of the bundle 
$\pi:{\mathbb C}{\times}{\mathcal E}^{(2)}
({\mathbb C}){\rightarrow}{\mathbb C}$, which is densely 
defined and multi-valued. In the next section, we investigate 
the solution of (\ref{parallel-equation}) more closely. 

\section{Bundle gerbes as non-cohomological notion}\label{3.0}

The bundle $\pi:{\mathbb C}{\times}{\mathcal E}({\mathbb C}){\rightarrow} 
{\mathbb C}$ 
together with the flat connection $\nabla$ gave us the notion of parallel
sections, where we extended this notion to be densely defined 
and multi-valued sections. 
This is in fact the notion of leafs of
the foliation given by the flat connection $\nabla$.
We now analyse the moduli space of the densely defined multi-valued 
parallel sections of the bundle 
$\pi:{\mathbb C}{\times}{\mathcal E}^{(2)}({\mathbb C}){\rightarrow}  
{\mathbb C}$ with respect to the connection $\nabla$. 
The moduli space has an unusual bundle structure, which we would 
call a {\it pile}. 
We analyse the evolution equation (\ref{parallel-equation}) for
parallel sections as a toy model of the phenomena for movable 
branching singularities.   

\subsection{Non-linear connections}\label{nonlinear}

First, consider a non-linear connection on the trivial bundle 
$\coprod_{\kappa\in\mathbb C}\mathbb C={\mathbb C}{\times}{\mathbb C}$ 
over ${\mathbb C}$ given by a holomorphic horizontal distribution 
\begin{equation}
  \label{eq:holdist}
H(\kappa; y)=
\{(t; y^2t); t\in{\mathbb C}\} \quad {\text{(independent of $\kappa$.)}}   
\end{equation}

The first equation of parallel translation \eqref{parallel-equation} 
is given by 
$\frac{dy}{d\kappa}=y^2$. Hence, parallel sections are given in general
\begin{equation}
 \label{eq:parallelsection}
(\kappa; y(\kappa)){=}(\kappa; \frac{1}{c{-}\kappa})
{=}(\kappa; \frac{c^{-1}}{1{-}c^{-1}\kappa}).  
\end{equation}
There are also singular solution $(\kappa; 0)$, corresponding
$c^{-1}=0$. Note that $(\kappa,-\frac{1}{\kappa})$ is not a singular solution.
For the consistency, we think that the singular point of the 
$(\kappa,0)$ section is $\infty$. 

Let ${\mathcal A}$ be the set of parallel sections 
including the singular solution $(\kappa,0)$.
Every $f{\in}{\mathcal A}$ has one singular point at 
a point $c \in S^2{=}{\mathbb C}\cup\{\infty\}$.     
The assignment of $f\in {\mathcal A}$ to its singular 
point $\sigma(f)=c$ gives a bijection 
$\sigma: {\mathcal A}{\to}S^2{=}{\mathbb C}\cup\{\infty\}$. 
Namely, ${\mathcal A}$ is parameterized by $S^2$ as    
\begin{equation}
 \label{eq:express}
\sigma(f)=c\Leftrightarrow f=(\kappa,\frac{1}{c{-}\kappa}),\quad  
\sigma(f)=\infty \Leftrightarrow f=(\kappa, 0)\in{\mathcal A}.  
\end{equation}
By this way, we give the topology on ${\mathcal A}$. 

\medskip
Let $T_{\kappa}^{\kappa'}(y)$ be the parallel displacement of
$(\kappa; y)$ along a curve from $\kappa$ to $\kappa'$. 
Since \eqref{eq:holdist} is independent of the base point
$\kappa$, $T_{\kappa}^{\kappa'}(y)$ is given by 
$$
T_{\kappa}^{\kappa'}(y)= \frac{y}{1{-}y(\kappa'{-}\kappa)}, \quad 
T_{\kappa}^{\kappa'}(\infty)=\frac{1}{\kappa{-}\kappa'}. 
$$
We easily see that 
$T_{\kappa}^{\kappa''}=T_{\kappa'}^{\kappa''}T_{\kappa}^{\kappa'},
\,\,T_{\kappa}^{\kappa}=I.$ 
Every $f{\in}{\mathcal A}$ satisfies 
$T_{\kappa}^{\kappa'}f(\kappa)=f(\kappa')$ where they are defined.

\subsubsection{Extension of nonlinear connection}\label{extended connection}

We now extend the non-linear connection $H$ defined by 
\eqref{eq:holdist} to the space 
${\mathbb C}{\times}{\mathbb C}^2$ by giving the holomorphic horizontal
distributions  
\begin{equation}
  \label{eq:holdist2}
{\tilde H}(\kappa; y, z)=\{(t; y^2t,-yt); t\in{\mathbb C}\} 
\quad {\text{(independent of $\kappa, z$)}}.  
\end{equation}
The parallel translation w.r.t. \eqref{eq:holdist2} is given by 
the following equations: 
\begin{equation}
  \label{eq:paratrans}
\frac{dy}{d\kappa}=y^2, \quad \frac{dz}{d\kappa}=-y.  
\end{equation}
For the equation \eqref{eq:paratrans}, multi-valued parallel sections
are given by two ways 
\begin{equation}
  \label{eq:parasec}
(\kappa, \frac{a}{1{-}a\kappa},z{+}\log(1{-}a\kappa), \quad 
(\kappa, \frac{1}{b{-}\kappa}, w{+}\log(\kappa{-}b)), \quad (a,b\in \C)  
\end{equation}
although they are infinitely many valued. 
Singular solutions $(\kappa; 0, z)$ 
are involved in the first expression.
Note that the set-to-set correspondence  
\begin{equation}
  \label{eq:settoset}
(a, z{+}{2\pi i}{\mathbb Z})\overset{\iota}{\Longleftrightarrow}
(b, w){=}(a^{-1},z{+}\log a{+}\pi i{+}{2\pi i}{\mathbb Z})  
\end{equation}
gives multi-valued parallel sections. However, because of 
the ambiguity of $\log a$, we can not make this correspondence 
a univalent correspondence. 

\medskip
Denote by ${\tilde{\mathcal A}}$ the set of all parallel sections  
written in the form \eqref{eq:parasec}. 
Denote by 
$\pi_3:{\tilde{\mathcal A}}{\to}{\mathcal A}$ be the forgetful 
mapping of the last component. This is surjective.
For every $v{\in}{\mathcal A}$ such that $\sigma(v){=}b{=}a^{-1}{\in}S^2$, we see 
$$
\pi_3^{-1}(v){=}
\{(\kappa, \frac{1}{b{-}\kappa}, w{+}\log(\kappa{-}b)); w{\in}{\mathbb C}\}{=} 
\{(\kappa,\frac{a}{1{-}a\kappa}, z{+}\log(1{-}a\kappa)); z{\in}{\mathbb C})\}. 
$$ 
Since there is one dimensional freedom of moving, 
$\pi_3^{-1}(v)$ should be parameterized by ${\mathbb C}$. 
However, there is no natural parameterization and there are 
many technical choices.  

\medskip

\subsubsection{Tangent spaces of ${\tilde{\mathcal A}}$}

For an element 
$f{=}(\kappa, \frac{a}{1{-}a\kappa}, z{+}\log(1{-}a\kappa))
{=}(\kappa,  \frac{1}{b{-}\kappa}, w{+}\log(\kappa{-}b))$, 
the {\it tangent space} $T_{f}{\tilde{\mathcal A}}$ of 
$\tilde{\mathcal A}$ at $f$ is  
$$
\begin{aligned}
T_{f}{\tilde{\mathcal A}}=&\big\{\frac{d}{ds}\Big|_{s=0}
(\frac{a(s)}{1{-}a(s)\kappa},z(s){+}\log(1{-}a(s)\kappa));(a(0),z(0)){=}(a,z)\big\}\\
=&\left\{(\frac{\dot{a}}{(1{-}a\kappa)^2}, \dot{z}{-}\frac{\dot{a}\kappa}{1{-}a\kappa}); 
\dot{a},\dot{z}{\in}{\mathbb C}\right\}
=\left\{(\frac{{-}\dot{b}}{(b{-}\kappa)^2}, \dot{w}{-}\frac{\dot{b}}{\kappa{-}b}); 
\dot{b},\dot{w}{\in}{\mathbb C}\right\}.
\end{aligned}
$$
Hence  
$$
\begin{bmatrix}
  \dot{b}\\\dot{w}
\end{bmatrix}{=}
\begin{bmatrix}
 - a^{-2}& 0\\
 a^{-1}& 1
\end{bmatrix}
\begin{bmatrix}
   \dot{a}\\\dot{z}
\end{bmatrix}{=}(d\iota)_{(a,z)}
\begin{bmatrix}
   \dot{a}\\\dot{z}
\end{bmatrix}. 
$$
Consider now a subspace $H_{f}$ of $T_{f}{\tilde{\mathcal A}}$
obtained by setting $\dot{z}=0$ in the definition of 
$T_{f}{\tilde{\mathcal A}}$.  
We regard $\{H_{f}; f\in\tilde{\mathcal A}\}$ as a horizontal
distribution on $\tilde{\mathcal A}$, where 
$\{H_{f}; f\in{\tilde{\mathcal A}}\}$ is defined without ambiguity.

\smallskip
The invariance for the vertical direction gives that 
$\{H_{f}; f\in{\tilde{\mathcal A}}\}$ is viewed as 
{\it infinitesimal trivialization} of 
$\pi_3:\tilde{\mathcal A}\to {\mathcal A}$.

\smallskip
The parallel translation $I_{\kappa}^{\kappa'}$ is given by 
\begin{equation}
  \begin{aligned}
I_{\kappa}^{\kappa'}(y,z)&{=}
(\frac{y}{1{-}y(\kappa'{-}\kappa)},z{+}\log(1{-}y(\kappa'{-}\kappa))),\\  
&{=}(\frac{1}{y^{-1}{-}\kappa'{+}\kappa},z{+}\log y{+}\log(y^{-1}{-}\kappa'{+}\kappa)),      
  \end{aligned}
\end{equation}
which is obtained by solving 
\eqref{eq:paratrans} under the initial data $(\kappa,y,z)$. 

By definition we see $I_{\kappa}^{\kappa}=I$, and 
$
I_{\kappa}^{\kappa''}=I_{\kappa'}^{\kappa''}I_{\kappa}^{\kappa'},
$ 
as a
set-to-set mapping 
Every $f\in {\tilde{\mathcal A}}$ satisfies 
$I_{\kappa}^{\kappa'}f(\kappa)=f(\kappa')$ where they are defined.

\begin{prop}
  \label{nongo}
The parallel displacement via 
horizontal distribution 
$\{H_f:f\in \tilde{\mathcal A}\}$ 
does not give a 
local trivialization of 
$\pi_3: {\widetilde{\mathcal A}}\to{\mathcal A}$.  
\end{prop}

\noindent
{\bf Proof}. For a point $f=(\kappa,\frac{a}{1{-}a\kappa})$ 
of ${\mathcal A}$, 
we consider a small neighborhood $V_a$ of $a$. 
Then 
${\tilde V}_a=\{\frac{a'}{1{-}a'\kappa}; a'\in V_a\}$ 
is a neighborhood of $f$ in ${\mathcal A}$. 
Consider the set 
$$
\pi_3^{-1}(V_a)=
\{(\kappa, \frac{a'}{1{-}a'\kappa},\log(1{-}a'\kappa)); a'\in V_a\},
$$ 
and the horizontal lift of the curve
$\frac{a'(s)}{1{-}a'(s)\kappa}$, $a'(s)=a{+}s(a'{-}a)$   
along the infinitesimal trivialization: Indeed this is to solve the equation  
$$
\frac{d}{ds}z(s)=\frac{(a'{-}a)\kappa}{1{-}a'(s)\kappa}, \quad z(0)\in\log(1{-}a\kappa).
$$
Hence 
$z(s)=\log(1{-}(a{+}s(a'{-}a))\kappa)$, and $z(1)=\log(1{-}a'\kappa)$. 
Thus the elimination of the ambiguity of $\log(1{-}a'\kappa)$ 
on the set $V_a$ is impossible, however small the neighborhood $V_a$
is. \qed

\medskip
Proposition\ref{nongo} shows that 
$\pi_3: {\widetilde{\mathcal A}}\to{\mathcal A}$ 
is not an affine bundle. In spite of this, one may say that the curvature of 
this connection vanishes.

\subsubsection{Affine bundle gerbe}

Although $\pi_3: {\widetilde{\mathcal A}}{\to}{\mathcal A}$ 
does not have a bundle structure, we can consider 
{\it local trivializations} by restricting the domain of $\kappa$.

\smallskip
\noindent
(a) Let $V_{\infty}{=}\{b; |b|{>}3\}\subset S^2$ be a neighborhood of $\infty$
First, we define a fiber preserving mapping $p_{_{\infty, D}}$ 
from the trivial bundle $\pi: V_{\infty}{\times}{\mathbb C}\to V_{\infty}$
into $\pi_3: {\widetilde{\mathcal A}}\to{\mathcal A}$ such that 
$\pi_3p_{_{\infty, D}}{=}\sigma^{-1}\pi$ by restricting the domain of 
$\kappa$ in a unit disk $D$: 
Consider $(\kappa, \frac{a}{1{-}a\kappa}, z{+}\log(1{-}a\kappa))$ for
$(\kappa, a^{-1}){\in}D{\times}V_{\infty}$.  
Since $|a\kappa|{<}1/3$, $\log(1{-}a\kappa)$ is defined as a univalent 
function $\log(1{-}a\kappa){=}\log|1{-}a\kappa|{+}i\theta$, 
${-}\pi{<}\theta{<}\pi$ on this domain by setting 
$1{-}a\kappa{=}|1{-}a\kappa|e^{i\theta}$, which will be denoted by 
$\log(1{-}a\kappa)_{D{\times}V_{\infty}}$.
We define   
\begin{equation}
  \label{eq:correp}
p_{_{\infty,D}}(b, z){=}(\kappa, \frac{a}{1{-}a\kappa},
z{+}\log(1{-}a\kappa)), \quad a^{-1}{=}b{\in}V_{\infty}, \quad
z{\in}{\mathbb C} 
\end{equation}
where $\log(1{-}a\kappa)$ in the right hand side is 
the analytic continuation of  
$\log(1{-}a\kappa)$ $=\log(1{-}a\kappa)_{D{\times}V_{\infty}}$.  

\smallskip
\noindent
\setlength{\unitlength}{.4mm}
\begin{picture}(140,100)
\put(10,5){\footnotesize{Fig.2}}
\put(50,50){\circle*{10}} 
\put(120,50){\circle*{10}}
\thinlines
\qbezier[100](20,50)(20,80)(50,80) 
\qbezier[100](20,50)(20,20)(50,20)
\qbezier[100](50,20)(80,20)(80,50)
\qbezier[100](80,50)(80,80)(50,80)
\thicklines
\qbezier[150](50,50)(85,67)(120,50)
\qbezier[150](50,50)(85,33)(120,50) 
\thinlines
\put(78,43){\circle{20}}
\put(78,57){\circle{20}}
\put(50,50){\circle{20}}
\put(72,30){\circle{20}}
\put(72,70){\circle{20}}
\put(60,22){\circle{20}}
\put(60,78){\circle{20}}
\put(58,35){\circle{20}}
\put(20,80){$V_{\infty}$}
\put(45,55){$D$}
\put(125,50){$D'$}
\put(30,42){$V_{b_0}$}
\put(85,33){$V_{b_1}$}
\put(84,62){$V_{b_{-1}}$}
\put(79,73){$V_{b_{-2}}$}
\put(79,20){$V_{b_{2}}$}
\put(60,12){$V_{b_{3}}$}
\put(107,58){$c_u$}
\put(107,38){$c_l$} 
\put(81,47){$\cdot$}
\put(80,46){$\cdot$}
\put(80,47){$\cdot$}
\put(80,48){$\cdot$}
\put(80,49){$\cdot$}
\put(81,48){$\cdot$} 
\put(82,47.5){$\cdot$}
\end{picture}\hfill
{\parbox[b]{.49\linewidth}{
(b) We take a simple covering of the domain $|z|\leq 3$ by unit
   disks $V_{b_{-k}},\dots, V_{b_{-1}},$ $V_{b_0},$ $V_{b_1},$
   $\dots,$ $V_{b_{\ell}}$ as in Fig.2, and fix a unit disk $D'$ apart from all $V_{b_i}$.      
We define a fiber preserving mapping 
$p_{_{V_{b_i},D'}}$ from the trivial bundle 
$\pi : V_{b_i}{\times}{\mathbb C}\to V_{b_i}$ into   
$\pi_3: {\widetilde{\mathcal A}}\to{\mathcal A}$ such that 
$\pi_3p_{_{V_{b_i},D'}}{=}\,\sigma^{-1}\pi$ by restricting the domain of $\kappa$ in 
a unit disk $D'$.} 

\smallskip
\noindent
We see that setting $\kappa{-}b{=}|\kappa{-}b|e^{i\theta}$, 
$\log(\kappa{-}b)$ is defined as a univalent function 
on the domain $D'{\times}V_{b_i}$ as 
$\log|\kappa{-}b|{+}i\theta$, $-\pi{<}\theta{<}\pi$, which is denoted by 
$\log(\kappa{-}b)_{D'{\times}V_{b_i}}$.

Consider $(\kappa,\frac{1}{b{-}\kappa}, w{+}\log(\kappa{-}b))$ for
$(\kappa, b){\in}D'{\times}V_{b_i}$. 
We define   
\begin{equation}
  \label{eq:correp2}
p_{_{V_{b_i},D'}}(b', w){=}(\kappa, \frac{1}{b'{-}\kappa},
w{+}\log(\kappa{-}b')), \quad (b', w){\in}V_b{\times}{\mathbb C}
\end{equation}
where $\log(\kappa{-}b)$ in the r.h.s. is the analytic
continuation of $\log(\kappa{-}b)_{V_{b_i}{\times}D'}$.  

\medskip
\noindent 
(c) Suppose $c{\in}V_{b_i}{\cap}V_{b_j}$ and 
$p_{V_{b_i},D'}(c, w){=}p_{V_{b_j}, D'}(c, w')$. Then we see that 
there exists uniquely $n(i,j){\in}{\mathbb Z}$ such that
$w'{=}w{+}2\pi in(i,j)$. For the above covering, we see 
$n(i,j){=}0$ for every pair $(i,j)$. 

\medskip
Let $c{\in}V_{b_i}{\cap}V_{\infty}$ and  
$p_{V_{b_i}, D'}(c, w){=}p_{{\infty}, D}(c, z)$. To fix the coordinate
transformation, we have to choose the identification 
of two sets of values $\log(\kappa{-}b)_{D'{\times}V_{b_i}}$ and 
$\log(1{-}a\kappa)_{D{\times}V_{\infty}}$. 
For $b_i$, except $b_1$, we identify these through the analytic
continuation along the (lower) curve $c_{l}$ joining $D$ and $D'$, but 
for $b_1$, we identify $\log(\kappa{-}b)_{D'{\times}V_{b_1}}$ and 
$\log(1{-}a\kappa)_{D{\times}V_{\infty}}$ through the analytic
continuation along the (upper) curve $c_{u}$ joining $D$ and $D'$.  

Therefore there is a positive integer $n(i,\infty)$ such that 
$w'{=}z{+}2\pi in(i,\infty)$ by the same argument. 
For the above covering we see in fact that if $n(1,\infty){=}\ell$ for $i{=}1$, then
$n(i,\infty){=}\ell{+}1$ for every $i{\not=}1$.  

These give coordinate transformations. However, 
the collection of these local trivializations do not glued together,  
for these do not satisfy the cocycle condition at the dotted small
triangle domain in Fig. 2. 

Thus, the collection of these local trivializations form only 
a {\bf bundle gerbe}. This may be denoted by $\coprod_{b{\in}S^2}{\mathbb C}_{b}$.
Thus we have a commutative diagram
\begin{equation}
 \label{eq:diag}
\begin{matrix}
\coprod_{b{\in}S^2}{\mathbb C}_{b}&
\overset{p_{(\bullet)}}{\longrightarrow}&\tilde{\mathcal A}\\
 \pi \downarrow & {}&  \pi_3\downarrow \\
{S^2}&\overset{\sigma}{\longleftarrow}&{\mathcal A}\\
\end{matrix}  
\end{equation}
One can consider various local trivializations of the 
bundle-like object of the left hand side.
$\coprod_{b{\in}S^2}{\mathbb C}_b$ is not an affine bundle, but 
an affine bundle gerbe with a holomorphic flat connection. 

Although $\tilde{\mathcal A}$ can not be a manifold, there are several 
manifold-like properties remain on the bundle-like object 
$\pi_3:{\tilde{\mathcal A}}{\to}{\mathcal A}$. 
On the affine bundle gerbe $\coprod_{b{\in}S^2}{\mathbb C}_{b}$, one 
can define a holomorphic flat connection via the diagram 
\eqref{eq:diag}. However, the geometric realization of a holomorphic 
parallel section is nothing but an element of $\widetilde{\mathcal{A}}$ 
given by \eqref{eq:parasec}.

\subsection{Geometric notions remained on $\tilde{\mathcal A}$}
Recall that the discordance (${sogo}$ in the term used in 
\cite{omori-proceedings}) of patching of three local coordinate 
neighbourhood occurs only on the small dotted triangle in Fig.2.  

In this section, we construct two examples which give almost 
the same phenomena as in the previous section for gluing local bundles.   
%
%
\subsubsection{Geometric quantization for a non integral 2-form}
Consider standard volume form $dV$ on $S^2$ with total volume $4\pi$. 
Let $\Omega$ be a closed smooth 2-form (current) on $S^2$ such that
$\int_{S^2}\Omega=4\pi\lambda$, but the support of $\Omega$
concentrates to a small disk neighborhood of the north pole $N$.   
Let $\{U_i\}_{i\in I}$ be a simple cover of $S^2$, on each $U_i$, 
$\Omega$ is written in the form $\Omega=d\omega_i$, and hence 
$\omega_{ij}=\omega_i{-}\omega_j$ on $U_{ij}=U_i\cap U_j$ is a closed 
1-form (current), which is written by $\omega_{ij}=df_{ij}$ on
$U_{ij}$ by using a smooth $0$-form (current) $f_{ij}$. 

Now we want to make a $U(1)$-vector bundle using $e^{\sqrt{-1}f_{ij}}$ as 
transition functions. However, since on $U_{ijk}=U_i\cap U_j\cap U_k$
we only have   
$$
e^{\sqrt{-1}f_{ij}}e^{\sqrt{-1}f_{jk}}e^{\sqrt{-1}f_{ki}}=
e^{\sqrt{-1}(f_{ij}{+}f_{jk}{+}f_{ki})},
$$  
$e^{\sqrt{-1}f_{ij}}$ can not be used as patching diffeomorphisms. 
In spite of these difficulties, we see that horizontal distributions 
defined by $\omega_i=0$ are glued together. 

Thus, we can define a linear connection on such a {\it broken} vector
bundle, which is indeed the notion of {\it bundle gerbes}. 
Since $\Omega=d\omega_i$, the curvature form of
this connection is given by $\Omega$. Remark that we can make a
parallel transform along any smooth curve $c(t)$ in $S^2$. 

Recall that the support of $\Omega$ is concentrated in a small
neighborhood $V_N$ of the north pole $N$. Therefore any closed 
curve in $S^2{-}V_N$ can shrink down to a point in  
$S^2{-}V_N$. In spite of this, the homotopy chasing of parallel 
transform does not succeed, because of the discordance (sogo) of the
patching diffeomorphisms.    

If $U_i$ does not intersect with $V_N$, then we have a product bundle 
$U_i\times{\mathbb C}$ with the trivial flat connection. 
Since $\omega_i=d\log e^{h_i}$, the integral submanifold of 
$\omega_i=0$ is given by $\log e^{h_i}$. 
This will be called the {\it pile}. 
Thus, even if the object is restricted to $S^2{-}V_N$, we have a 
{\it non-trivial} bundle gerbe which is apparently not classified 
by cohomologies. 

We note that this gives also a concrete example of local line bundles 
treated by \cite{melrose} over a manifold. 

\subsubsection{A simple example}
The simplest example of objects we propose in this paper is given
by the Hopf-fibering $S^3\,{\stackrel{S^1}{\rightarrow}}\,S^2$.
Viewing $S^3{=}\coprod_{q\in S^2}S_q^1$ (disjoint union),
we consider the $\ell$-covering $\tilde{S}_q^1$ of each fiber 
${S}_q^1$, and denote by $\tilde{S}^3$ the disjoint union 
$\coprod_{q\in S^2}\tilde{S}_q^1$.
We are able to define local trivializations of
$\tilde{S}^3|_{U_i}\cong U_i{\times}\tilde{S}^1$ naturally through the
trivializations ${S}^3|_{U_i}$ given on a simple open covering 
$\{U_i\}_{i\in \Gamma}$ of $S^2$. This structure permits
us to treat $\tilde{S}^3$ as a local Lie group, and hence it
looks like a topological space. On the other hand, we have a projection
$$
\pi \,\,:\,\,\tilde{S}^3{=}\!\coprod_{q\in S^2}\tilde{S}_q^1\rightarrow S^3{=}
\!\coprod_{q\in S^2}{S}_q^1
$$
as the union of fiberwise projections, as if it were a non-trivial $\ell$-
covering. However $\tilde{S}^3$ cannot be a manifold, since $S^3$ is
simply connected. In particular, the {\it points} of $\tilde{S}^3$ 
should be regarded as {\it $\ell$-valued elements}. 

\par
We now consider a 1-parameter subgroup $S^1$ of $S^3$ and the
inverse image $\pi^{-1}(S^1)$. Since all points of ${\widetilde S}^3$
are ``$\ell$-valued'', this simply looks like a combined object of 
$S^1\times {\mathbb Z}_{\ell}$ and the $\ell$ covering group, i.e. in some
restricted region, this object can be regarded as a point set by
several ways. In such a region, the ambiguity is caused simply by 
the reason that two pictures of point sets are mixed up.

\subsubsection{Conceptual difficulties beyond ordinary mathematics}
Let $P_c$ be the parallel translation along a closed curve. Let
$c_s(t)$ be a family of closed curves. Suppose 
$c_s(0)=c_s(1)=p$ and $c_1(t)=p$. We see that there is 
$(p;v)$ such that $P_{c_s}(p;v)\not=v$. Therefore there must be somewhere 
a singular point for the homotopy chasing, caused by the discordance.    
However the position of singular point can not be specified. 

Even though the parallel transformation is defined for every 
fixed curve, these parallel transformations are in general set-to-set
mappings when one parameter family of closed curves are considered. 

Thus, we have some conceptual difficulty that may be explained as follows: 
a parallel transformation along a curve has a definite meaning, 
but when we think this in a family of curves, then we have to think 
{\it suddenly} this a set-to-set mapping. Recall here the 
``Schr{\"o}dinger's cat''.  

Recall that such a strange phenomena caused in $\tilde{\mathcal A}$ by   
movable branching singularities. 
In \S\ref{nonlinear}, we considered a nonlinear connection on the 
trivial bundle $S^2{\times}{\mathbb C}$, and an extended connection 
to treat the amplitude of the star exponential functions of the quadratic 
form. 

\section{Broken associative product and extensions}

In this section we give an example where such a fuzzy phenomenon play 
a crucial role to define a concrete algebraic structure. 
We consider the product bundle $\coprod_{\kappa\in{S^2}}{\mathbb C}$, and we define  
at each fiber another associative product, which is {\it broken} in the sense 
that they are not necessarily defined for all pairs $(a,b)$. 

\subsection{Associative product combined with Cayley transform}

First of all we give such a product on the fiber at $\kappa{=}0$. 
Let $S^2$ be a 2-sphere identified with $\C\cup\{\infty\}$. 
Consider the Cayley transform 
$C_0: S^2{\to}S^2$, $C_0(X)=\frac{1{-}X}{1{+}X}$, and define the
product by  
\begin{equation}
  \label{eq:tanh}
a{\bullet}_0b=\frac{a{+}b}{1{+}ab}\sim C_0^{-1}(C_0(a)C_0(b)).  
\end{equation}
Here $\sim$ means the equality under the algebraic calculation, which is 
an algebraic procedure through the calculations such as follows: 
$$\frac{1-\frac{1-a}{1+a} \cdot \frac{1-b}{1+b}}
{1+\frac{1-a}{1+a} \cdot \frac{1-b}{1+b}}\sim
\frac{(1{+}a)(1{+}b)-(1{-}a)(1{-}b)}{(1{+}a)(1{+}b)+(1{-}a)(1{-}b)}
= \frac{a+b}{1+ab}.
$$
This is defined for every pair $(a,b)$ such that $ab\not=-1$, but commutative 
and associative whenever they are defined.  
Remark also that 
\begin{equation}
  \label{eq:reminv}
a{\bullet}_0b=\frac{a{+}b}{1{+}ab}\sim\frac{a^{-1}{+}b^{-1}}{1{+}(ab)^{-1}}
=a^{-1}{\bullet}_0b^{-1}.  
\end{equation}
Hence we set 
$\infty{\bullet}_0b=b^{-1},\,\,
\infty{\bullet}_0\infty{=}0 $ in particular.

One can extend this broken product to pairs 
$(a;g){\in}{\C}{\times}{\C}$ as follows: 
$$
(a:g){\bullet}_0(b:g')=(a{\bullet}_0b: gg'(1{+}ab)). 
$$
This is an associative product. The associativity gives 
$$ 
(1{+}bc)(1{+}a\frac{b{+}c}{1{+}bc})=(1{+}\frac{a{+}b}{1{+}ab}c)(1{+}ab). 
$$
It is worthwhile to write this identity in the logarithmic form
\begin{equation}
  \label{eq:logcocycl}
\log(1{+}bc){+}\log(1{+}a\frac{b{+}c}{1{+}bc})=
\log(1{+}\frac{a{+}b}{1{+}ab}c){+}\log(1{+}ab), \quad
{\text{mod}}\,2\pi{i}{\mathbb Z} 
\end{equation}
although the logarithmic form uses infinitely many valued functions. 
If one denote $C(a, b)=\log(1{+}ab)$, then the above identity means 
the Hochschild's 2-cocycle condition: 
$$
C(b,c){-}C(a{\bullet}_0b,c){+}C(a, b{\bullet}_0c){-}C(a,b)=0, \quad 
{\text{mod}}\,2\pi{i}{\mathbb Z}.
$$

We extend the product as follows:
\begin{equation}
\label{lnprod}
(a:g){\bullet}_{ln}(b:g')=(a{\bullet}_0b: g{+}g'{+}\log(1{+}ab)).  
\end{equation}
This is associative as a set-to-set mappings.  
Thus, we extends \eqref{lnprod} for all pairs $a,b\in S^2$
except $(0,\infty)$ by setting  
$$
(\infty:g){\bullet}_{ln}(b:g')=(b^{-1}:g{+}g')=(b:g){\bullet}_{ln}(\infty,g').
$$ 

\smallskip
Next we define a family of products defined on each fiber at $\kappa$.
To define such a product, we use the twisted Cayley transform  
defined by $C_{\kappa}\sim C_0T^0_{\kappa}$, where $T^0_{\kappa}$ is
given in the equality \eqref{intertwiner}. The result is   
\begin{equation}
  \label{eq:twstcly}
C_{\kappa}(y)=\frac{1{-}y(1{-}\kappa)}{1{+}y(1{+}\kappa)},  
\end{equation}
and we  define 
\begin{equation}
  \label{eq:prod2}
a{\bullet}_{\kappa}b=\frac{a{+}b{+}2ab\kappa}{1{+}ab(1{-}\kappa^2)}
\sim C_{\kappa}^{-1}(C_{\kappa}(a)C_{\kappa}(b)).
\end{equation}
The point is that the singular set of product changes by the algebraic
trick of calculations.  
$a{\bullet}_{\kappa}b$ is defined for every pair $(a,b)$ such that 
$ab(1{-}\kappa^2){\not=}{-}1$. Another word, for an arbitrary pair 
$(a,b)\in {\mathbb C}^2$, the product 
$a{\bullet}_{\kappa}b$ is defined for some $\kappa$ moving in 
an open dense domain. 

For parallel sections given in \eqref{eq:parallelsection} we see  
\begin{equation}
  \label{eq:paraprod}
\frac{a}{1{-}a\kappa}{\bullet}_{\kappa}\frac{b}{1{-}b\kappa}=
\frac{a{+}b}{1{-}(a{+}b)\kappa{+}ab}.  
\end{equation}
In particular,  
$$
{\kappa}^{-1}{\bullet}_{\kappa}{\kappa}^{-1}=0, \quad 
{\kappa}^{-1}{\bullet}_{\kappa}\frac{1}{b^{-1}{-}\kappa}=\frac{1}{b{-}\kappa}.
$$
For simplicity, we denote by $f(\kappa)$ the section $f$ of the bundle 
$\pi S^2\times \C \rightarrow S^2$. 

\begin{prop}
 \label{q-func}
For parallel sections $f(\kappa), g(\kappa)$ defined on open subsets, 
the product $f(\kappa){\bullet}_{\kappa}g(\kappa)$ is also a parallel
section where they are defined.    
\end{prop}

\subsubsection{Extension of the product}
Using \eqref{lnprod}, 
one can extend the product $a{\bullet}_{\kappa}b$ by the formula 
$$
(a; g){\bullet}_{\kappa}(b; g')\sim 
I_0^{\kappa}\big((I_{\kappa}^{0}(a; g)){\bullet}_{ln}(I_{\kappa}^{0}(a; g'))\big).
$$
Indeed, we see how the algebraic trick works: 
\begin{equation}
  \label{eq:extprod}
  \begin{aligned}
(a:g){\bullet}_{\kappa}(b:g')=&(a{\bullet}_{\kappa}b:g{+}g'{+}\log(1{+}ab(1{-}\kappa^2))\\
=&\Big(\frac{a{+}b{+}2ab\kappa}{1{+}ab(1{-}\kappa^2)}:
 g{+}g'{+}\log(1{+}ab(1{-}\kappa^2))\Big).     
  \end{aligned}
\end{equation}

\begin{prop}
  \label{extprod2000}
The extended product $(a:g){\bullet}_{\kappa}(b:g')$ is defined 
with the ambiguity of $2\pi i{\mathbb Z}$. However, ${\bullet}_{\kappa}$ product is 
associative whenever these are defined.
\end{prop}
The point of such a fiberwise product is the following: 
\begin{prop}
 \label{q-func22}
For parallel sections $f(\kappa), g(\kappa)$ defined on open subsets, 
the product $f(\kappa){\bullet}_{\kappa}g(\kappa)$ is also a parallel
section where they are defined.    
\end{prop}

\noindent
{\bf Proof}. We have only to prove 
$I_{\kappa}^{\kappa'}(f{\bullet}_{\kappa}h)=
I_{\kappa}^{\kappa'}(f){\bullet}_{\kappa'}I_{\kappa}^{\kappa'}(h). $  

For $f=(\frac{a}{1{-}a\kappa},\log(1{-}a\kappa)),\, 
h=(\frac{b}{1{-}b\kappa},\log(1{-}b\kappa))$, we see that 
$$
\begin{aligned}
f{\bullet}_{\kappa}h=&
\Big(\frac{a{+}b}{1{-}(a{+}b)\kappa{+}ab},  
\log\big((1{-}a\kappa)(1{-}b\kappa)
(1{+}\frac{a}{1{-}a\kappa}\frac{b}{1{-}b\kappa}(1{-}\kappa^2))\big)\Big)\\
=&\Big(\frac{a{+}b}{1{-}(a{+}b)\kappa{+}ab},\log(1{-}(a{+}b)\kappa{+}ab)\Big).   
\end{aligned}
$$
It is easily seen that 
$I_{\kappa}^{\kappa'}(f{\bullet}_{\kappa}h)=
\Big(\frac{a{+}b}{1{-}(a{+}b)\kappa'{+}ab}, \log(1{-}(a{+}b)\kappa'{+}ab)\Big). $
 \qed

\section{Notion of $q$-number functions}

Using Propositions\,\ref{q-func}, \ref{q-func22}, we define 
a multiplicative structure on the sets $\mathcal A$ and 
$\tilde{\mathcal A}$ of all parallel sections. 
A notion of 
$q$-number functions which 
describe a quantum observables was introduced in \cite{AW}, where 
our notion of the parallel sections 
is stimulated by this idea. 
From this point of view, we may employ the notation 
$:f:_\kappa$ for a section $f$  
of the bundle $\pi:S^2\times \C\rightarrow S^2$. 

We think $f{\in}{\mathcal A}$ is an element viewing $\kappa$ as an indeterminate. 
For every $f, g{\in}{\mathcal A}$ excluding the pair  
$(f,g)=(\frac{1}{1{-}\kappa},\frac{-1}{1{+}\kappa})$ we define an element 
$f{\bullet}g{\in}{\mathcal A}$ by  
\begin{equation}
  \label{eq:prodprod}
{:}f{\bullet}g{:}_{\kappa}=f({\kappa}){\bullet}_{\kappa}g(\kappa).  
\end{equation}

Some of product formula on $\mathcal A$ are given as follows:
$$
0{\bullet}f=f,\quad \frac{-1}{\kappa}{\bullet}\frac{-1}{\kappa}=0, 
\quad \frac{1}{1{-}\kappa}{\bullet}f=\frac{1}{1{-}\kappa}, \quad 
\frac{-1}{1{+}\kappa}{\bullet}f=\frac{-1}{1{+}\kappa}, 
$$
where $0$ stands for the singular solution $(\kappa, 0)$.
These say that $\frac{\pm 1}{1{\mp}\kappa}$ looks like $0$ or
$\infty$. Hence ${\mathcal A}$ is viewed naturally as the Riemann
sphere with standard multiplicative structure such that 
$a0=0$, $a\infty=\infty$, but $0\infty$ does not defined. 
By the definition of $\bullet_{\kappa}$, we have 
$C_{\kappa}(f{\bullet}_{\kappa}g)=C_{\kappa}(f)C_{\kappa}(g)$. 

Hence the correspondence is given by the family of 
twisted Cayley transform 
$\coprod_{\kappa\in {\mathbb C}}C_{\kappa}:
{\mathcal A}{\to}{\mathbb C}\,{\cup}\{\infty\}$. 
We view $\mathcal A$ as a topological space through the 
identification $\coprod_{\kappa{\in}{\mathbb C}}C_{\kappa}$.  

The table of correspondence is as follows: 
\medskip
\begin{center}\label{lnprod5}
\begin{tabular}{c|ccccc|cc}
${\mathcal A}$ &0& $\frac{-1}{\kappa}$&$\frac{1}{1{-}\kappa}$
&$\frac{-1}{1{+}\kappa}$&$\frac{a}{1{-}a\kappa}$&
$\frac{1{-}a}{1{-}\kappa{+}a(1{+}\kappa)}$&$f(\kappa)$\\\hline
${\rm{Im}}{C_{\kappa}}$&$1$&$-1$&$0$&$\infty$&$\frac{1{-}a}{1{+}a}$&
$a$&$\frac{1{-}f(\kappa)(1{-}\kappa)}{1{+}f(\kappa)(1{+}\kappa)}$\\
$\sigma$;$\kappa=$&$\infty$&$0$&$1$&$-1$&$\frac{1}{a}$&$\frac{1{+}a}{1{-}a}$& --\\
\end{tabular}
\end{center}
Note that 
$$
C_{\kappa}^{-1}(a)=
\frac{1{-}a}{1{-}\kappa{+}a(1{+}\kappa)}\sim 
\frac{\frac{1{-}a}{1{+}a}}{1{-}\frac{1{-}a}{1{+}a}\kappa}
\sim T_0^{\kappa}C_0^{-1}(a)
$$
is a parallel section, and 
$\frac{1{-}f(\kappa)(1{-}\kappa)}{1{+}f(\kappa)(1{+}\kappa)}$ is free
from $\kappa$ for every parallel section.

\subsection{Product on $\tilde{\mathcal A}$}

Let ${\tilde{\mathcal A}}$ be the space of all parallel sections 
given in \eqref{eq:parasec}, and consider the product 
${\bullet}$ on $\tilde{\mathcal A}$ is given by 
the product formula \eqref{eq:extprod}. For $f, f'\in\tilde{\mathcal A}$,
we set $f=(\kappa,y(\kappa),z(\kappa))$, $f'=(\kappa,y'(\kappa),z'(\kappa))$.
$f{\bullet}g$ is defined as a parallel section defined on the open dense 
domain such that $y(\kappa), y'(\kappa)\not=\infty$. 

Remark that 
$$
(\kappa,\frac{a}{1{-}a\kappa}, \log(1{-}a\kappa)){\bullet}
(\kappa, \frac{-1}{1{+}\kappa},\log(1{+}\kappa))=
(\kappa,\frac{-1}{1{+}\kappa},\log(1{-}a){+}\log(1{+}\kappa)),
$$
$$
(\kappa,\frac{a}{1{-}a\kappa},\log(1{-}a\kappa)){\bullet}
(\kappa,\frac{1}{1{-}\kappa},\log(1{-}\kappa))=
(\kappa,\frac{1}{1{-}\kappa},\log(1{+}a){+}\log(1{-}\kappa)).
$$
Although $\frac{\pm1}{1{\mp}\kappa}$ plays a role of $0$ or $\infty$, the
second component depends on $a$. 

For simplicity, we denote in particular  
\begin{equation}
  \label{eq:vacvacvac}
\varpi_{c}=(\kappa; \frac{1}{1{-}\kappa}, c{+}\log\frac{1}{2}(1{-}\kappa)),
\quad
\bar\varpi_{c}=(\kappa; \frac{-1}{1{+}\kappa}, c{+}\log\frac{1}{2}(1{+}\kappa)).  
\end{equation}
It is easy to see that 
$$
\varpi_{c}{\bullet}\varpi_{c'}=\varpi_{c{+}c'},\quad
\bar\varpi_{c}{\bullet}\bar\varpi_{c'}=\bar\varpi_{c{+}c'}, 
$$
but $\varpi_{c}{\bullet}\bar\varpi_{c'}$ diverges.

Let 
${\tilde{\mathcal A}}_{\times}$ be the subset of 
${\tilde{\mathcal A}}$ excluding parallel sections  
$(\kappa; \frac{\pm 1}{1{\mp}\kappa}, c{+}\log(1{\mp}\kappa))$.
We denote also 
$$
{\tilde{\mathcal A}}_{0}={\tilde{\mathcal A}}_{\times}\cup\{\varpi_{c}\},\quad
{\tilde{\mathcal A}}_{\infty}={\tilde{\mathcal A}}_{\times}\cup\{\bar\varpi_{c}\}.
$$

\begin{prop}
\label{paraprod22}
${\widetilde{\mathcal A}}$ is closed under the extended product
$\bullet_{\kappa}$, whenever these are defined. 
In particular, 
${\widetilde{\mathcal A}}_{\times}$, ${\widetilde{\mathcal A}}_{0}$
${\widetilde{\mathcal A}}_{\infty}$ are closed respectively under the 
$\bullet$-product.   
\end{prop}

\subsection{Infinitesimal left-action}

Note that the singular solution 
${\bf 1}{=}(\kappa,0,0)\in{\tilde{\mathcal A}}$ is the multiplicative
identity. A neighborhood of ${\bf 1}$ is given as 
$(\kappa, \frac{a}{1{-}a\kappa},g{+}\log(1{-}a\kappa))$ by taking $(a,g)$
in a small neighborhood of $0$. Setting $g=0$,
we denote $f_a=(\kappa, \frac{a}{1{-}a\kappa},\log(1{-}a\kappa))$. 
For a parallel section 
$h=(\kappa, y(\kappa), z(\kappa))\in{\mathcal A}$, 
the product $f_a{\bullet}h$ is given  
$$
f_a{\bullet}h=
(\kappa, 1{+}y(\kappa)(1{+}\kappa){-}y(\kappa)^2(1{-}\kappa^2)),  
$$

Consider the infinitesimal action 
$$
\frac{d}{ds}\Big|_{s=0}(as ,0){\bullet}_{\kappa}(y,z)=
\frac{d}{ds}\Big|_{s=0}(a(1{+}2y\kappa{-}y^2(1{-}\kappa^2)), ay(1{-}\kappa^2)).
$$
Define for every fixed $\kappa$ the invariant distribution 
$$
{\tilde{L}}_{\kappa}(y,z)=\{(a((1{+}y\kappa)^2{-}y^2), ay(1{-}\kappa^2)); a
\in {\mathbb C}\}. 
$$

By Proposition\ref{paraprod22}, we have 
$dI_0^{\kappa}{\tilde{L}}_{0}I_{\kappa}^0(y,z)={\tilde{L}}_{\kappa}$.

\subsection{Exponential mapping}
The equation of integral curves of the invariant distribution 
$ {\tilde{L}}_{\kappa} $
through the identity $(0,0)$ is  
$$
\frac{d}{dt}(y(t), z(t))=(a((1{+}y(t)\kappa)^2{-}y(t)^2),
ay(t)(1{-}\kappa^2)), \quad (y(0),z(0))=(0,0).
$$
For the case $\kappa=0$, $a=1$, we have 
$(y(t), z(t))= (\tanh t, \log\cosh t).$

\medskip
We define ${\rm{Exp}}_{\bullet}: {\mathbb C}\to {\mathcal A}_{\times}$ by 
the family of ${\rm{Exp}}_{\kappa}$: 
$$
{\rm{Exp}}_{\bullet_{\kappa}}t =T_0^{\kappa}(\tanh t)=
\frac{\sinh t}{\cosh t{-}(\sinh t)\kappa},
$$
\begin{equation}
  \label{eq:expkappa}
{\rm{Exp}}_{\bullet}t =(\kappa; T_0^{\kappa}(\tanh t))=
\Big(\kappa; \frac{\sinh t}{\cosh t{-}(\sinh t)\kappa}\Big).  
\end{equation}

For a fixed $t$, ${\rm{Exp}}_{\bullet}t$ is a parallel section with the exponential law  
$$
{\rm{Exp}}_{\bullet}s{\bullet}{\rm{Exp}}_{\bullet}t={\rm{Exp}}_{\bullet}(s{+}t),
 \quad{\text{and}}\quad
{\rm{Exp}}_{\bullet}(s{+}2\pi i)={\rm{Exp}}_{\bullet}s.
$$ 

\medskip
For the extended product, let ${\widetilde{{\rm{Exp}}}}_0\,t=(\tanh t; \log\cosh t)$,  
and let 
$$
{\widetilde{{\rm{Exp}}}}_{\kappa}\,t= I_0^{\kappa}{\widetilde{{\rm{Exp}}}}_0\,t=
\Big(\frac{\sinh t}{\cosh t{-}(\sinh t)\kappa}, \log(\cosh t{-}(\sinh t)\kappa\Big).
$$ 
Although ${\widetilde{{\rm{Exp}}}}_{\kappa}$ is not defined for all 
$t\in {\mathbb C}$, viewing $\kappa$ as an indeterminate permits us to 
define the exponential mapping  
${\widetilde{{\rm{Exp}}}}_{\bullet}: {\mathbb C}\to{\widetilde{\mathcal A}}_{\times}$ by 
\begin{equation}
  \label{eq:expgen}
{\widetilde{{\rm{Exp}}}}_{\bullet}t=
(\kappa; 
\frac{\sinh t}{\cosh t{-}(\sinh t)\kappa}, \log(\cosh t{-}(\sinh t)\kappa)).  
\end{equation}
This is a parallel section with the exponential law  
$$
{\widetilde{{\rm{Exp}}}}_{\bullet}s{\bullet}
{\widetilde{{\rm{Exp}}}}_{\bullet}\,t={\widetilde{{\rm{Exp}}}}_{\bullet}(s{+}t).
$$
Remark here that ${\widetilde{{\rm{Exp}}}}_{\bullet}s$ is infinitely many valued. 
We have ${\widetilde{{\rm{Exp}}}}_{\bullet}(s{+}2\pi i)={\widetilde{{\rm{Exp}}}}_{\bullet}s$,
but these equalities hold by set-to-set correspondence. 

Consider the formula  
\begin{equation}
  \label{eq:intlog}
\log(\cosh t{-}(\sinh t)\kappa)=
\int_0^{t}
\frac{\sinh s{-}(\cosh s)\kappa}{\cosh s{-}(\sinh s)\kappa}ds.  
\end{equation}
The multi-valuedness of the left hand side is caused by the 
integral depends on the homology class of the route of integration. 
Thus, for a fixed $\kappa$, the left hand side must be considered 
as a univalent function on the homological universal covering space of 
${\mathbb C}{-}\{\cosh s{-}(\sinh s)\kappa=0\}$ (the path space
factored by the group of all homologically trivial loops).  

Minding this, we see the following   
\begin{prop}
  \label{logmul}
If $\kappa$ is treated as an indeterminate, then 
$2\pi i$ periodicity does not appear in the integral 
$$
\int_0^{t}
\frac{\sinh s{-}(\cosh s)\kappa}{\cosh s{-}(\sinh s)\kappa}ds.
$$
By this, ${\widetilde{{\rm{Exp}}}}_{\bullet}$ must be viewed as an injective homomorphism.   
\end{prop}

\noindent
{\bf Proof}\,\,If the periodicity appears, then 
$\int_a^{a{+}2\pi i}
\frac{\sinh s{-}(\cosh s)\kappa}{\cosh s{-}(\sinh s)\kappa}ds$ 
must be the one of $2\pi n$. Hence the derivative by $\kappa$ must 
vanish. 
Consider now the following quantity: 
$$
\frac{d}{d\kappa}\int_a^{a{+}2\pi i}
\frac{\sinh s{-}(\cosh s)\kappa}{\cosh s{-}(\sinh s)\kappa}ds= 
-\int_a^{a{+}2\pi i}\frac{1}{(\cosh s{-}(\sinh s)\kappa)^2}ds.
$$
This does not vanish if we regard $\kappa$ as an indeterminate.  \qed

\medskip
We see also that for every $\alpha\in{\mathbb C}$
$$
{\widetilde{{\rm{Exp}}}}^{(\alpha)}_{\bullet}s=
\Big(\kappa;\frac{\sinh s}{\cosh t{-}(\sinh s)\kappa}, 
\log{e^{\alpha s}}(\cosh s{-}(\sinh s)\kappa)\Big)
$$
satisfies the exponential law.
Using this formula, it is easily seen that for $t\in{\mathbb R}$,
$$
{:}\varpi_{0}{:}_{\kappa}= 
\lim_{t\to\infty}
\Big(\frac{\sinh t}{\cosh t{-}(\sinh t)\kappa}, 
\log{e^{-t}}(\cosh t{-}(\sinh t)\kappa)\Big), 
$$
$$
{:}\bar\varpi_{0}{:}_{\kappa}= 
\lim_{t\to -\infty}
\Big(\frac{\sinh t}{\cosh t{-}(\sinh t)\kappa},
\log{e^{t}}(\cosh t{-}(\sinh t)\kappa)\Big).
$$

\bigskip
$\tilde{\mathcal A}_{\times}$ is a strange object. This is a 
group-like object. The forgetful mapping is a homomorphism 
onto ${\mathcal A}_{\times}\cong{\mathbb C}_{\times}$.
There is an injective exponential mapping 
${\widetilde{{\rm{Exp}}}}_{\bullet}: {\mathbb C}\to\tilde{\mathcal A}_{\times}$.
In spite of this, one can not treat $\widetilde{\mathcal A}_{\times}$
as a manifold. 


\end{document}